\documentstyle[12pt]{article}

\newcommand{\R}{\mbox{$I\!\!R$}}

\begin{document}

\title{On regular algebraic surfaces of $\R^3$ with constant mean curvature}
\author{J. Lucas M. Barbosa and Manfredo P. do Carmo \thanks{%
Both authors are partially supported by CNPq, Brazil}}
\date{}
\maketitle
\begin{abstract}
We consider regular surfaces $M$ that are given as the zeros of a polynomial
function $p:\R^3\rightarrow \R$, where the gradient of $p$ vanishes nowhere.
We assume that $M$ has non-zero mean curvature and prove that there exist
only two examples of such surfaces, namely the sphere and the circular cylinder.
\end{abstract}

\section{Introduction}
An algebraic set in $\R^3$ will be the set
$$
M = \{ (x,y,z) \in \R^3;\; p(x,y,z)=0\}
$$
of zeros of a polynomial function $p:\R^3\rightarrow R$. An algebraic set is regular
if the gradient vector $\nabla p = (p_x, p_y, p_z)$ vanishes nowhere in $M$; here
$p_x$, $p_y$ and $p_z$ denotes the derivative of $p$ with respect to $x$, $y$ or $z$ respectively.

The condiction of regularity is essential in our case. It allows to parametrize the set $M$
locally by differentiable functions $x(u,v)$, $y(u,v)$, $z(u,v)$ (not necessarily polynomials),
so that $M$ becomes a regular surface in the sense of differential geometry ( see \cite{dc} chapter 2
section 2.2, in particular Proposition 2); here $(u,v)$ are coordinates in an open set of $\R^2$.

Since $M$ is a closed set in $\R^3$, it is a complete surface. In addition, being a regular surface,
it is properly embedded, i. e., the limit set of $M$ (if any) does not belong to $M$ (cf. \cite{Wh},
chapter IV, A.1 p. 113). In particular, regular algebraic surfaces are locally graphs over their
tangent planes,

From now on, $M$ will denote a regular algebraic surface in $\R^3$. Due to the regularity condition,
one can define on $M$ the basic objects of Differential Geometry of surfaces and pose some
differential-algebraic questions within this algebraic category.

For instance, in the last 60 years (namely after the seminal work \cite{Ho1} of Heinz Hopf in 1951),
many questions have been worked out on differentiable. surfaces of non-zero constant mean curvature
$H$. See also \cite{Ho2}.

In our case, we have two examples of algebraic regular surfaces that have non-zero constant
mean curvature, namely,
\begin{enumerate}
\item [(1)] Spheres, $(x-x_0)^2 + (y - y_0)^2 + (z - z_0)^2 = r^2$ with center $(x_0, y_0, z_0) \in \R^3$ and
radius $r = 1/H$.
\item[(2)] Circular right cylinders, $(x-x_0)^2 + (y - y_0)^2 = r^2$, whose basis is a circle in the plane $xy$
with center $(x_0, y_0)$ and whose axis is a straight line passing through the center and parallel to the
$z$ axis.
\end{enumerate}
A first natural question is: Are there further examples?

\vspace{3mm}

The first time we heard about this question was in a preprint of Oscar Perdomo (recently published in \cite{pe})
where he proves that for polynomials of degree three there are no such surfaces.

\vspace{3mm}

\noindent In this note, we prove the following general result:

\vspace{3mm}

\noindent {\bf Theorem:} {\em Let $M$ be a regular algebraic surface in $\R^3$. Assume that it has constant mean curvature $H \ne 0$. Then
$M$ is a sphere or a right circular cylinder.}

\vspace{2mm}

\noindent {\bf Acknowledgements.} We want to thank Fernando Cod\'a Marques for a crucial observation, Oscar Perdomo for having written \cite{pe}
and Karl Otto St\"ohr for his help in our first attempts to solve the problem.

\section{Preliminaries}

We first observe that, in the compact case, this theorem follows immediately from Alexandrov's well-known result:
{\it An embedded compact surface in $\R^3$ with constant mean curvature is isometric to a sphere.}

The second observation is that the total curvature of an algebraic surface is finite. This was first proved
by Osserman \cite{Os} in the case that the surface is an immersion parametrized by polynomials in two variables.

Here we give a proof of the finiteness of the total curvature for a (implicitly defined) regular algebraic surface $M$.

\vspace{3mm}

\noindent {\bf Proposition $1$:} {\em \label{2prop1}
Let $M$ be a regular algebraic surface in $\R^3$. Assume that it is complete non-compact. Then its total curvature
is finite. More precisely,
$$
\int_M |K|\, dM \leq 4\pi C(d),
$$
where $K$ is the Gaussian curvature of $M$ and $C(d)$ is a constant that depends only on the degree $d$ of the polynomial
defining $M$.}

\vspace{3mm}

\noindent {\bf Proof.} Let $g:M \rightarrow S^2(1)$ be the Gauss map of $M$. It is well known the Gauss curvature
$K = \det (dg)$, where $dg$ is the differential of the map $g$. Of course the result is true if $M$ is either a plane
or a circular cylinder. So, we ruled out these two cases from our proof. Let $M^*$ be the set of points in $M$ where
$K\neq 0$. Then, restricted to $M^*$, $g$ is a local diffeomorphism; that is, given $q\in g(M^*)$ and $m_{\alpha} \in
\{ g^{-1}(q)\}$, $\alpha$ belonging to a set of indices $A$, there exist neighborhoods $\cal U$ of $q$ and ${\cal V}_{\alpha}$ of $m_{\alpha}$
such that, for each $\alpha$, $g$ maps $V_{\alpha}$ diffeomorphically onto $\cal U$. In fact, $g$ restrict to $M^*$ is a covering
of $N(M^*)$ without ramification points. Since
$$
\int_M |K| dM = \int_{M^*} |K|dM
$$
the theorem is proved if we show that the mentioned covering has only a finite number of leaves.

Now, fix a plane $P$ passing through the origin of $\R^3$ and in $P$ fix an orthonormal basis $\{e_1, e_2\}$. Then, there exists a point $m\in M$ such that, up to translations, $P = T_m(M)$ if and only if
$$
\frac{\nabla p}{|\nabla p|}(m)\; \bot \; P
$$
what is equivalent to the system of equations
\begin{equation} \left\{ \begin{array} {cc} \label{basic}
 \langle \frac{\nabla p}{|\nabla p|}(m), e_1 \rangle = 0& \\
\langle \frac{\nabla p}{|\nabla p|}(m), e_2 \rangle = 0&
\end{array}\right.
\end{equation}
Set $e_1 = \sum a_iU_i$ and $e_2 = \sum b_iU_i$ where $U_1 = (1,0,0)$, $U_2 =(0,1,0)$, $U_3 = (0,0,1)$ is the canonical basis of $\R^3$. Then, since $|\nabla p| \ne 0$, then (\ref{basic}) takes the form
\begin{equation}\left\{\begin{array} {cc}\label{basic1} 
p_x a_1 + p_y a_2 + p_z a_3 = 0& \\
p_x b_1 + p_y b_2 + p_z b_3 = 0&
\end{array} \right.
\end{equation}
where the $a_i, b_i$, $i= 1,2,3$ are real numbers.

The equations in (\ref{basic1}) describe algebraic surfaces, $\Sigma_1$ and
$\Sigma_2$, which correspond to the coefficients $a_i$ and $b_i$
respectively, and whose degrees are $\le (d-1)$, where $d$ is the degree of $p$.

The surfaces $\Sigma_1$ and $\Sigma_2$, together with the original surface $M$
determines points $m \in M$ as follows:

$\Sigma_j$ interset $M$, $j=1,2$, in a curve $C_j$. If the intersection $C_1\cap C_2$ contains a point $m\in M^*$, since $K(m)\ne 0$, such intersection is unique in a
neighborhood of $m$. By Bezout's theorem the total number of intersections is bounded above by $(d-1)^2$, as we wished. This proves the proposition.

\vspace{3mm}

Since $\int_M |K|\,dM$ is finite, it follows from a theorem of Huber \cite{Huber} that the surface $M$ of Propositon $1$ is finitely connected, i.e.,
it is a compact surface with a finite number of ends.

The proof of our Theorem uses in a crucial way the structure theory for embedded, complete finitely connected surfaces with non-zero constant mean curvature
developed by Korevaar, Kusner and Solomon in \cite{KKS} after some preliminary work by Meeks \cite{Meeks}. The statement that we need from these papers is as follows:

\vspace{3mm}

\noindent {\bf Theorem A:} (\cite{Meeks} and \cite{KKS}.) {\em Let $M$ be a complete, non-compact, properly ebedded surface in $\R^3$ with non-zero constante
mean curvature. Assume that $M$ is finitely connected. Then, the ends of $M$ are cylindrically bounded. Furthermore, for each end $E$ of $M$, there exists
a Delaunay surface $\Sigma \subset \R^3$ such that $E$ and $\Sigma$ can be expressed as cylindrical graphs $\rho_E$ and $\rho_{\Sigma}$
so that, near infinity, $|\rho_E - \rho_{\Sigma}|<C e^{-\lambda x }$ where $C\ge 0$ and $\lambda >0$ are constants.}

\vspace{3mm}

\noindent {\bf Remark:} The first assertion in Theorem A comes from \cite{Meeks}. The final assertion is from \cite{KKS}, theorem 5.18.

\section{Proof of the Theorem}

We can assume that $M$ is complete and non-compact; otherwise it is a sphere. Thus, by Proposition 1, $M$ has finite total curvature, and hence, by Huber's theorem, $M$ is compact with finitely many ends. By Theorem A, each end $E$ of $M$ converges exponentially to a Delaunay surface $\Sigma$.
Since $M$ is embedded, the Delaunay surface $\Sigma$ to which an end $E$ converges has to be an onduloid or a right circular cylinder.

We first claim that the Delaunay surface $\Sigma$ towards which $E$ converges is actually a cylinder.

Suppose it is not. By a rigid motion, we can assume that the axis of $\Sigma $ is parallel to the $y$ axis and meets the $z$ axis.
Then, there is a value $z_0$
of $z$ such that the line $y \rightarrow (0, y, z_0)$ intersects $\Sigma $ infinitely often. Since $E$ approaches $\Sigma$ at infinity,
the algebraic equation $p(0,y,z_0) = 0$ has infinitely many solutions. This is impossible. So $\Sigma$ is a cylinder as we claimed.

We claim now that $E$ contains a open set of the cylinder $\Sigma$.

To see this, we take a rigid motion so that one of the straight lines of the cylinder $\Sigma$ agrees with the coordinate $y$-axis.
Thus, one of the intersection curves of $E$ with the plane $x=0$ is a curve $\beta$ that converges to the $y$-axis. If $y$ is large
enough, $\beta$ is given by
$$
\beta(y) = (0, y, z(y)),
$$
where $z(y)$ is a function that satisfies
$$
\lim_{y \rightarrow \infty} z(y) = 0
$$
Since the curve $\beta$ belongs to the end $E$, we have
\begin{equation} \label{eq5} \nonumber
p(0,y,z(y)) =0
\end{equation}
Observe that the polynomial $p$ can be written as
$$
p(x,y,z) = a_nz^n + a_{n-1}z^{n-1} + \ldots + a_1z + a_0
$$
where $a_k = a_k(x,y)$ is a polynomial in $x$ and $y$ of degree $\le n$.
By Theorem A, we have that
$$
\lim_{y\rightarrow \infty} z(y) = \lim_{y\rightarrow \infty} Ce^{-\lambda y} = 0
$$
By a known result in Calculus, we have, for any integer $k$,
$$
\lim_{y\rightarrow \infty} y^k e^{-\lambda y} = 0
$$
for any integer $k$.

Thus, by computing the limit in the equation $p(0,y,z(y)) = 0$ as $y\rightarrow \infty$ along the curve $\beta$,
we obtain that $a_0$ does not depend on $Y$, and $a_0=0$. This means that, for any $y$, the equation
$p(0,y,z)= 0$ has $z=0$ as a root, i.e., the straight line $y\rightarrow (0,y,0)$ is contained in $E$.

The above argument applies to an arbitrary straight line of $\Sigma$. It follows that an open set in $E$ is a cylinder. This proves our claim

Thus, there exists an open set $U$ in $M$ with the property that the Gaussian curvature $K$ vanishes in $M$. Since $M$ is analytic, $K$
vanishes identically in $M$. It is then well known (see e.g. \cite{Ma}) that $M$ is a cylinder. Since $H$ is constant, this is a
circular cylinder. This proves the Theorem.

\vspace{3mm}

\noindent {\bf Remark:} A crucial point in the proof is that the convergence in \cite{KKS} is exponential. It allows us to prove that not only
an arbitrary line in the cylinder $\Sigma$ converges to $E$ but that actually it is contained in $E$.

\section{Final Remarks}

{\bf The case $H=0$.} There are many algebraic minimal surfaces in $\R^3$ (see p. 161 of the English translation of Nitsche's
book \cite{Ni}). However, the examples we are most familiar with, namely, the Enneper surface and the Hennenberg surface, are not embedded; thus
they are not regular algebraic surfaces.

In fact it is simple to prove the following proposition

\vspace{2mm}

\noindent {\bf Proposition:} {\em There are no regular algebraic minimal surfaces in $\R^3$ except the plane. }

\vspace{2mm}

\noindent {\bf Proof:} Let $M$ be an algebraic surface in $\R^3$. As we have seen, such surface it is finitely connected, i.e., it is a compact surface with a finite number of ends. We also know that $M$ is properly embedded.

Let $E$ be one of its ends. Parametrically $E$ can be described by a map
$x:D - \{O\}\rightarrow \R^3$, where $D$ is an open disk of $\R^2$ centered at the origin and $O$ is the origin.

 We may assume, after a rotation if necessary, that the Gauss map, which extends to $O$ (see Osserman \cite{Os1}), takes on the value $(0, 0, 1)$ at $O$. The two simplest examples of such ends are the plane and (either end of) the catenoid.

Now we use a result proved by R. Schoen \cite{Sc}. He showed that such an end is the graph of the function $x_3$ defined over the $(x_1, x_2)$-plane and

When $a \ne 0$ the end is of catenoid type. When $a=0$ the end of the planar type. In fact, if $a\ne 0$ the function $x_3$ will be asymptotic to the graph of the function $\log \rho$; if $a=0$ it will be asymptotic to the graph of a constant function (equal to $\beta$).

Let's assume that $E$ is of the catenoid type.
Consider the curve $\alpha$ intersection of the $E$ with the plane $x_2=0$ in the region $x_1>0$.
Since $M$ is given by the equation $p(x_1,x_2,x_3)=0$, then the curve $\alpha$
is algebraic, given by $p(x_1,0,x_3) =0$. This curve must be asymptotic
to the graph of the function $x_3 = a\log x_1$. But this is impossible. Hence, $M$ can not have end of the catenoid type.

Thus, all the ends of $M$ are of the planar type. But they are in finite number. Since $M$ is embedded, the planes asymptotic to $M$ must be parallel. It follows that there are two parallel planes such that $M$ is contained in the region bounded by them. It follows by the halfspace theorem for minimal surfaces \cite{HM} that $M$ must be a plane.

\vspace{2mm}

\noindent {\bf Hypersurfaces in $\R^{n+1}$, $n\ge 3$.} In this case we consider the zeros of a polynomial function
$p(x_0, x_1,\ldots, x_n)$, $n\ge 3$, with $\nabla p \ne 0$ everywhere, and call it {\it regular algebraic hypersurfaces} $M^n$ of
$\R^{n+1}$. Similar to the case $n=2$, the only compact example of such hypersurfaces are spheres. This follows immediately from
Alexandrov theorem. So, we are left to consider the complete non-compact case. A generalized cylinder $C^k$ in $\R^{n+1}$ is a
product $B^k \times \R^{n-k}$, where the basis $B^k \subset \R^{k+1}\subset \R^{n+1}$ is a hypersurface of $\R^{k+1}$ and the
product is embedded in $\R^{n+1}$ in the canonical way, i.e.,  $B^k\times \R^{n-k} \subset \R^{k+1}\times \R^{n-k}$. It is
easily checked that when $B$ is a $k$-sphere, $C^k$ has nonzero constant mean curvature. The following lemma is again a
consequence of Alexandrov's theorem.

\vspace{3mm}

\noindent {\bf Lemma.} {\em Let $C^k$ be an algebraic regular generalized cylinder in $\R^{n+1}$ whose basis $B$ is a compact hypersurface.
If $C^k$ has constant mean curvature then the bases $B^k$ is a $k$-sphere.}

\vspace{3mm}

We do not know any further examples of a regular algebraic hypersurfaces in $\R^{n+1}$, $n>2$, with nonzero constant mean curvature. We can
ask a question similar to the one we answered for $n=2$. The possible extension of our proof, however, needs new ideas. Although the total
Gauss-Kronecker curvature is again finite, there is no Huber theorem for $n>2$, and the proof of the structure theorem of \cite{KKS}
does not work por hypersurfes in $\R^{n+1}$, $n>2$.

\vspace{3cm}

\noindent Manfredo Perdig\~ao do Carmo \\
Instituto Nacional de Matem\'atica Pura e Aplicada - IMPA \\
Estrada Dona Castorina 110 \\
22460-320 Rio de Janeiro - RJ \\
Brazil \\
manfredo@impa.br

\vspace{1cm}

\noindent Jo\~ao Lucas Marques Barbosa \\
Rua Carolina Sucupira 723 Ap 2002 \\
60140-120 Fortaleza - Ce \\
Brazil \\
joaolucasbarbosa@terra.com.br

\end{document}